\documentclass[a4paper,12pt]{article}

\usepackage[top=30truemm,bottom=25truemm,left=25truemm,right=25truemm]{geometry}
\usepackage[pdftex]{graphicx}
\usepackage{ulem}
\renewcommand{\title}[1]{
	\vspace*{4mm}
	\begin{center}
	\textbf{\Large #1}
	\end{center}
	\smallskip}

\renewcommand{\author}[1]{
	\vspace*{0mm}
	\begin{center}
	#1
	\end{center}
	\smallskip}
\renewcommand{\abstract}[1]{
	\begin{center}
	\parbox{13cm}{\small {\sc Abstract.} #1}
	\end{center}
	\smallskip}

\renewcommand{\thanks}[1]{{
	\renewcommand{\thefootnote}{\fnsymbol{footnote}}
	\vspace*{-5mm}
	\footnote[0]{#1}}}
\newcommand{\address}[1]{\bigskip{\small\noindent #1 \par}}
\newcommand{\email}[1]{{\small\noindent\textit{Email address}: \texttt{#1} \par}}

\usepackage{titlesec} 
\titleformat*{\section}{\large\bfseries}
\titleformat*{\subsection}{\bfseries}
\titleformat*{\subsubsection}{}
\titlespacing{\section}{0pt}{*3}{*1.5}
\titlespacing{\subsection}{0pt}{*2}{*1}
\titlespacing{\subsubsection}{0pt}{*2}{*1}

 \usepackage{amscd, amsmath, amssymb, amsthm}
 \usepackage[initials,short-journals,non-sorted-cites]{amsrefs}
 \usepackage{tikz}
\usetikzlibrary{
  knots,
  hobby,
  decorations.pathreplacing,
  shapes.geometric,
  calc
}

\usepackage{graphicx}
\usepackage{color}

 \theoremstyle{plain}
 \newtheorem{thm}{Theorem}[section]
 \newtheorem{prop}[thm]{Proposition}
 
 \newtheorem{lem}[thm]{Lemma}
 \newtheorem{question}[thm]{Question}

 \theoremstyle{definition}
 \newtheorem{defn}[thm]{Definition}
 \newtheorem{example}[thm]{Example}
 \newtheorem{rem}[thm]{Remark}

\begin{document}

\tikzset{
  knot diagram/every strand/.append style={
    ultra thick,
    black 
  },
  show curve controls/.style={
    postaction=decorate,
    decoration={show path construction,
      curveto code={
        \draw [blue, dashed]
        (\tikzinputsegmentfirst)--(\tikzinputsegmentsupporta)
        node [at end, draw, solid, red, inner sep=2pt]{};
        \draw [blue, dashed]
        (\tikzinputsegmentsupportb)--(\tikzinputsegmentlast)
        node [at start, draw, solid, red, inner sep=2pt]{}
        node [at end, fill, blue, ellipse, inner sep=2pt]{}
        ;
      }
    }
  },
  show curve endpoints/.style={
    postaction=decorate,
    decoration={show path construction,
      curveto code={
        \node [fill, blue, ellipse, inner sep=2pt] at (\tikzinputsegmentlast) {}
        ;
      }
    }
  }
}
\tikzset{->-/.style 2 args={
    postaction={decorate},
    decoration={markings, mark=at position #1 with {\arrow[thick, #2]{>}}}
    },
    ->-/.default={0.5}{}
}

\tikzset{-<-/.style 2 args={
    postaction={decorate},
    decoration={markings, mark=at position #1 with {\arrow[thick, #2]{<}}}
    },
    -<-/.default={0.5}{}
}
\title{A construction of multiple group racks}

\author{Katsunori Arai}

 \abstract{
  A multiple group rack is a rack which is a disjoint union of groups equipped with
  a binary operation satisfying some conditions.
  It is used to define invariants of spatial surfaces, i.e., oriented compact surfaces with boundaries embedded in the $3$-sphere $S^{3}$.
  A $G$-family of racks is a set with a family of binary operations indexed by the elements of a group $G$.
  \textcolor{black}{There are two known methods for constructing multiple group racks.
  One is via a $G$-family of racks. The resulting multiple group rack is called the associated multiple group rack of the $G$-family of racks.
  The other is by taking an abelian extension of a multiple group rack.}
  In this paper,
  we introduce a new method for constructing multiple group racks by using a $G$-family of racks and a normal subgroup $N$ of $G$.
  \textcolor{black}{
  We show that this construction yields multiple group racks 
  that are neither the associated multiple group racks of any $G$-family of racks nor
  their abelian extensions when the right conjugation action of $G$ on $N$ is nontrivial.
  }
  As an application, 
  we present a pair of spatial surfaces 
  that \textcolor{black}{cannot be distinguished by invariants derived from the associated multiple group racks of any $G$-family of racks,
  yet} can be distinguished using invariants obtained from a multiple group rack introduced in this paper.
 }




\section{Introduction}
It is well known that every oriented link in the $3$-sphere $S^{3}$ bounds a compact oriented surface.
Such surfaces for a given oriented link are called \textit{Seifert surfaces}.
Classifying Seifert surfaces for a given link up to ambient isotopy in the link complement is an important problem in knot theory.
A \textit{spatial surface} is a compact surface embedded in $S^{3}$.
Throughout this paper,
we assume that (1) a spatial surface is oriented and 
(2) each connected component of a spatial surface is neither a disk nor a surface without boundaries.
Under the assumptions,
every spatial surface is a Seifert surface for its boundaries.
Two spatial surfaces are said to be 
\textit{equivalent} if they are ambiently isotopic in $S^{3}$.
When two spatial surfaces with the same boundaries are not equivalent,
then they are not ambiently isotopic in the complement of their boundaries.
Consequently,
classifying spatial surfaces contributes to the classification of Seifert surfaces.
A \textit{spatial trivalent graph} is a finite trivalent graph embedded in $S^{3}$.
In this paper, we allow graphs to have multiple edges, loops, and \textit{circle components}, that is, edges without vertices.
Diagrams of spatial trivalent graphs are defined analogously to knot diagrams.
In \cite{Matsuzaki2021},
a presentation of spatial surfaces using spatial trivalent graph diagrams was introduced,
and 
a Reidemeister-type theorem for spatial surfaces was established.

A \textit{quandle} \cite{Joyce1982,Matveev1982} is an algebraic system whose axioms correspond to Reidemeister moves in knot theory.
Given a quandle,
a \textit{coloring} of an oriented link diagram is a map from the set of all arcs in the diagram to the quandle, 
satisfying a specific condition at each crossing.
It is well known that the number of colorings 
is an invariant of oriented links.
A \textit{multiple conjugation quandle} \cite{Ishii2015'} is a quandle 
whose underlying set is a disjoint union of groups, and 
it is used to construct invariants of \textit{handlebody-knots},
i.e., handlebodies embedded in $S^{3}$.
A \textit{rack} \cite{Fenn-Rourke1992} is a generalization of a quandle and 
\textcolor{black}{its axioms correspond} to two of the three Reidemeister moves.
A \textit{multiple group rack} \cite{Ishii-Matsuzaki-Murao2020} is a rack defined on a disjoint union of groups equipped with a binary operation 
that satisfies axioms corresponding to the Reidemeister moves for spatial surfaces.
Just as the axioms of a rack are obtained by removing one of the axioms of a quandle,
the axioms of a multiple group rack are obtained by relaxing one of the axioms of a multiple conjugation quandle.
We can construct invariants of spatial surfaces by using multiple group racks, 
such as the number of colorings of spatial surface diagrams by given multiple group racks.
Therefore,
these invariants 
can be used to distinguish spatial surfaces. 
Several examples of distinguishing spatial surfaces using invariants derived from multiple group racks
can be found in \cite{Ishii-Matsuzaki-Murao2020,Matsuzaki-Murao2023}.

\textcolor{black}{A \textit{$G$-family of racks} \cite{Ishii-Iwakiri-Jang-Oshiro2013} is a set equipped with binary operations on the set, indexed by the elements of a group $G$, satisfying some conditions.
There are two known methods for constructing multiple group racks.
One is via a $G$-family of racks.
The resulting multiple group rack is called the \textit{associated multiple group rack} of the $G$-family of racks \cite{Ishii2015',Ishii-Matsuzaki-Murao2020}.
The other is by taking an \textit{abelian extension} \cite{Carter-Ishii-Saito-Tanaka2017,Matsuzaki-Murao2023} of a given multiple group rack.
In the present paper,
we introduce a new method for constructing multiple group racks by using $G$-families of racks and a normal subgroup $N$ of the group $G$.
We show that this construction yields multiple group racks that are neither the associated multiple group racks of any $G$-family of racks 
nor their abelian extensions when the right conjugation action of $G$ on $N$ is nontrivial.}
Moreover,
we provide a pair of spatial surfaces that satisfy the following conditions:
(1) their boundaries are ambiently isotopic oriented links,
(2) their regular neighborhoods in $S^{3}$ are ambiently isotopic handlebody-knots, 
(3) they cannot be distinguished by the number of colorings using the associated multiple group racks of any $G$-family of racks, and
and (4) they can be distinguished by the number of colorings by multiple group racks in the present paper.

In section~\ref{Sec:Spatial_surfaces},
we \textcolor{black}{recall} how to present spatial surfaces via diagrams of spatial trivalent graphs
and check the Reidemeister-type theorem for spatial surfaces.
Additionally,
we discuss \textcolor{black}{the notion of Y-orientations} of diagrams,
which is necessary for defining colorings of spatial surface diagrams by using multiple group racks.
In section~\ref{Sec:MGR}, 
we \textcolor{black}{recall} the definition of \textcolor{black}{a multiple group rack} and describe 
the construction of the associated multiple group racks of $G$-families of racks.
In section~\ref{Sec:MGRCols}, 
we recall the definition of a coloring of Y-oriented diagrams by multiple group racks.
Moreover,
we see a key property of 
colorings of Y-oriented diagrams by the associated multiple group rack of a $G$-family of racks.
In section~\ref{Sec:Main_results},
we present a new construction for multiple group racks by using a $G$-family of racks and a normal subgroup of the group $G$ (Theorem~\ref{Thm:Main_result}).
We also prove that the multiple group racks obtained by this method satisfy a certain property (Theorem~\ref{Thm:Main_result2}).
As an application,
we provide a pair of spatial surfaces 
that can be distinguished using an invariant obtained from a multiple group rack constructed in the way of Theorem~\ref{Thm:Main_result},
even though their boundaries and regular neighborhoods in $S^{3}$ are, respectively, ambiently isotopic.
Furthermore, 
the coloring number of the associated multiple group racks of any $G$-family of racks fails to distinguish them.

\section{Diagrams of Spatial surfaces}{\label{Sec:Spatial_surfaces}}

A \textit{spatial surface} is a compact surface embedded in $S^{3} = \mathbb{R}^{3} \sqcup \left\{\infty\right\}$.
Throughout this paper, we assume the following conditions: 
(1) a spatial surface is oriented, and (2) each connected component of a spatial surface is neither a disk nor a surface without boundaries.
Two spatial surfaces are said to be \textit{equivalent} if they are ambiently isotopic in $S^{3}$. 
For two spatial surfaces $F_{1}$ and $F_{2}$, 
we use the notation $F_{1} \cong F_{2}$ when $F_{1}$ and $F_{2}$ are equivalent.
A \textit{spatial trivalent graph} is a finite trivalent graph embedded in $S^{3}$.
In this paper, we allow graphs to have multiple edges, loops, and \textit{circle components}, that is, edges without vertices.
Diagrams of spatial trivalent graphs are defined analogously to knot diagrams.
An \textit{edge} of a diagram $D$ of a spatial trivalent graph $G$ is a sub-diagram of $D$ presenting an edge of $G$.
In particular,
an edge presenting a circle component of $G$ is also called a \textit{circle component} of $D$.

\begin{figure}[h]
  \centering
  \includegraphics{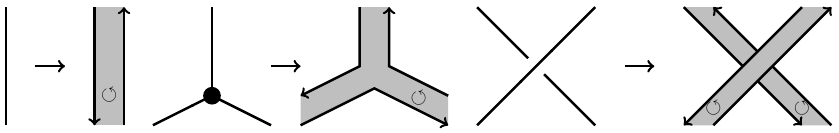}
  \caption{A construction of a spatial surface from a diagram of a spatial trivalent graph}
  {\label{Fig:Construction_of_spatial_surfaces}}
\end{figure}

\begin{figure}[h]
  \centering
  \includegraphics{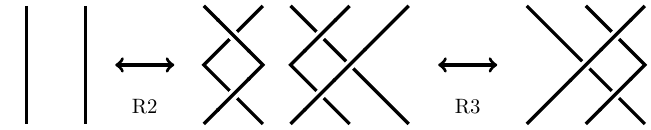}
  
  \vspace*{5mm}
  \includegraphics{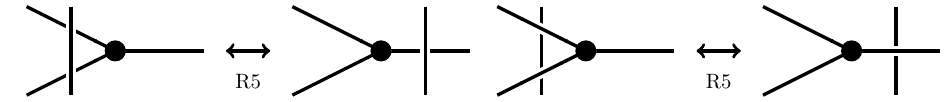}

  \vspace*{5mm}
  \includegraphics{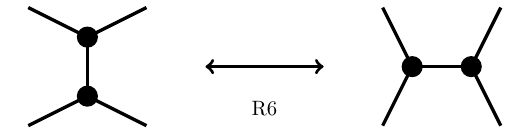}
    \caption{Reidemeister moves for diagrams of spatial surfaces}{\label{Fig:R-moves}}
\end{figure}

Let $D$ be a diagram of a spatial trivalent graph.  
Consider the spatial surface $F(D)$ obtained from $D$ as illustrated in Fig.~\ref{Fig:Construction_of_spatial_surfaces}. 
More precisely, take a regular neighborhood $N(D)$ of $D$ in $\mathbb{R}^{2}$ and 
replace it locally around each crossing with two bands in $\mathbb{R}^{3}$, as shown in the rightmost part of Fig.~\ref{Fig:Construction_of_spatial_surfaces}. 
This process yields a compact surface embedded in $\mathbb{R}^{3}$.
Assign an orientation to the surface induced by the orientation of $\mathbb{R}^{2}$.
By considering this oriented surface in $S^{3} = \mathbb{R}^{3} \sqcup \left\{\infty\right\}$,
we obtain the spatial surface $F(D)$, 
referred to as the \textit{spatial surface obtained from $D$}.
For any spatial surface $F$,
there exists a spatial trivalent graph diagram $D$ such that $F \cong F(D)$, \cite{Ishii-Matsuzaki-Murao2020,Matsuzaki2021}.
A \textit{diagram} of a spatial surface $F$ is a spatial trivalent graph diagram $D$ such that $F(D)$ is 
equivalent to $F$.
As in the case of knots, 
a Reidemeister-type theorem holds for spatial surfaces.
\thm[\cite{Matsuzaki2021}]{\label{Thm:R-moves}}{
	Two spatial surfaces are equivalent if and only if
	their diagrams are related by a finite sequence of $\mathrm{R}2$, $\mathrm{R}3$, $\mathrm{R}5$ and $\mathrm{R}6$ moves, as depicted in Fig.~\ref{Fig:R-moves}, and isotopies in $S^2 = \mathbb{R}^{2} \cup \{ \infty \}$.   
}
\upshape

\begin{figure}[h]
  \centering
  \includegraphics{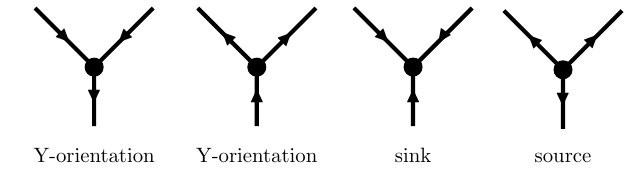}
  \caption{All orientations around trivalent vertices}{\label{Fig:All_orientation}}
\end{figure}

\begin{figure}[h]
  \centering
\includegraphics{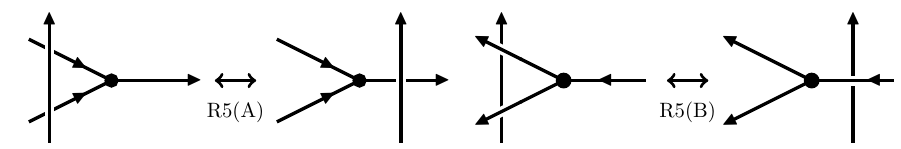}

\includegraphics{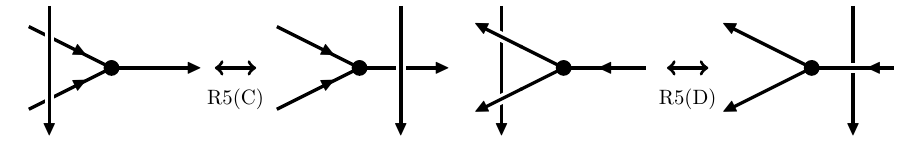}

\includegraphics{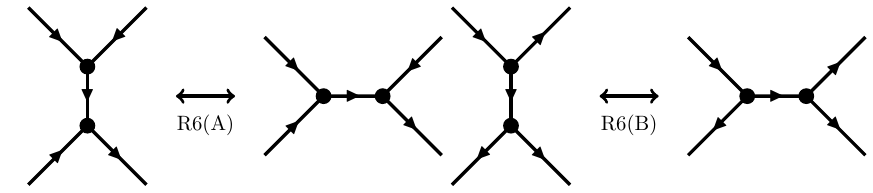}

\includegraphics{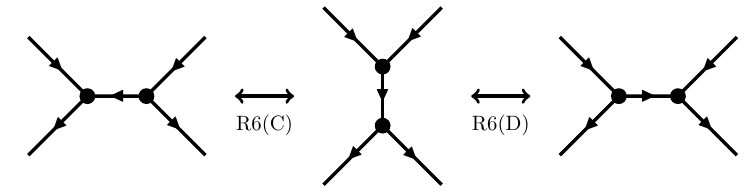}

\includegraphics{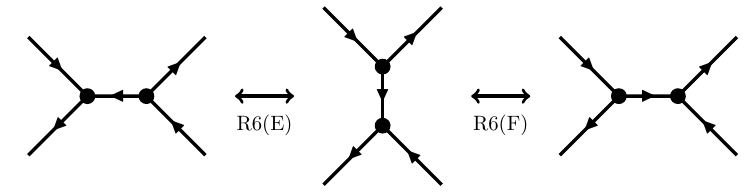}
  \caption{Y-oriented R$5$ and R$6$ moves}{\label{Fig:Y-oriented_R-moves}}
 \end{figure}

Let $D$ be a diagram of a spatial surface.
A \textit{Y-orientation} of $D$ is an assignment of orientations to all edges of $D$
such that no vertices are sinks or sources, as shown in 
Fig.~\ref{Fig:All_orientation}.
It is known that any diagram 
admits at least one Y-orientation, \cite{Ishii2015,Lebed2015}.
A \textit{Y-oriented diagram} of a spatial surface is a diagram equipped with a Y-orientation.
\textit{Y-oriented $\mathrm{R}5$ and $\mathrm{R}6$ moves} are local moves on Y-oriented diagrams whose underlying moves are 
$\mathrm{R}5$ and $\mathrm{R}6$ moves, as illustrated in Fig.~\ref{Fig:Y-oriented_R-moves}.
\textit{Y-oriented Reidemeister moves} for Y-oriented diagrams of spatial surfaces are local moves
generated by oriented $\mathrm{R}2$ and $\mathrm{R}3$ moves and Y-oriented $\mathrm{R}5$ and $\mathrm{R}6$ moves.

Two Y-oriented diagrams $D$ and $D'$ with the same underlying diagram are related
by a finite sequence of Y-oriented Reidemeister moves, reversing some circle components, and isotopies in $S^{2}$, \cite{Ishii2015,Matsuzaki-Murao2023}.
Thus, a Reidemeister-type theorem for Y-oriented diagrams follows from Theorem~\ref{Thm:R-moves}.
\thm[\cite{Matsuzaki-Murao2023}]{\label{thm Y-oriented R-moves}}{
  Two spatial surfaces are equivalent if and only if 
  their Y-oriented diagrams are related by a finite sequence of Y-oriented Reidemeister moves,
  reversing some circle components, and 
  isotopies in $S^2$.
}
\upshape

\section{Multiple group racks}{\label{Sec:MGR}}

\defn[\cite{Fenn-Rourke1992,Joyce1982,Matveev1982}]{\label{Def:rack}}{
  A \textit{rack} $X$ is a set with a binary operation $\ast : X \times X \to X$ satisfying the following conditions.
  \begin{itemize}
    \item[(i)] For any $y \in X$, the map $S_{y} : X \times X \to X$, defined by $x \mapsto x \ast y$, is bijective.
    \item[(ii)] For any $x, y, z \in X$, $(x \ast y) \ast z = (x \ast z) \ast (y \ast z)$.
  \end{itemize}
  A rack $X$ is called a \textit{quandle} if the binary operation $\ast$ satisfies $x \ast x = x$ for any $x \in X$.
}



\defn[\cite{Ishii-Matsuzaki-Murao2020}]{\label{Def:MGR}}{
  A \textit{multiple group rack} $X = \bigsqcup_{\lambda \in \Lambda} G_{\lambda}$ is a disjoint union of groups $G_{\lambda}$ ($\lambda \in \Lambda$)
  with a binary operation $\ast : X \times X \to X$ satisfying the following conditions.
  \begin{itemize}
    \item[(i)] For any $x \in X$, $\lambda \in \Lambda$, and $a, b \in G_{\lambda}$,
    $x \ast (ab) = (x \ast a) \ast b$ and $x \ast e_{\lambda} = x$,
    where $e_{\lambda}$ is the identity element of the group $G_{\lambda}$.
    \item[(ii)] For any $x, y, z \in X$,
    $(x \ast y) \ast z = (x \ast z) \ast (y \ast z)$.
    \item[(iii)] For any $\lambda \in \Lambda$ and $x \in X$,
    there is an element $\mu \in \Lambda$ such that
      for any $a, b \in G_{\lambda}$, $a \ast x, b \ast x \in G_{\mu}$ and $(ab) \ast x = (a \ast x) (b \ast x)$.
  \end{itemize}
  \textcolor{black}{A \textit{multiple conjugation quandle} \cite{Ishii2015'} is a multiple group rack $X = \bigsqcup_{\lambda \in \Lambda} G_{\lambda}$ such that}
  $X$ satisfies $a \ast b = b^{-1} ab$ for any $a, b \in G_{\lambda}$.
}



\defn[\cite{Ishii-Iwakiri-Jang-Oshiro2013}]{\label{Def:G-family}}{
  Let $G$ be a group with the identity element $e$.
  A \textit{$G$-family of racks} is a set $X$ with a family of binary operations $\left\{\ast^{g} : X \times X \to X\right\}_{g \in G}$ satisfying the following conditions.
  \begin{itemize}
    \item[(i)] For any $x, y \in X$ and $g, h \in G$,
    $x \ast^{gh} y = (x \ast^{g} y) \ast^{h} y$ and $x \ast^{e} y = x$.
    \item[(ii)] For any $x, y, z \in X$ and $g, h \in G$,
    $(x \ast^{g} y) \ast^{h} z = (x \ast^{h} z) \ast^{h^{-1} gh} (y \ast^{h} z)$.
  \end{itemize}
}

\example{\label{Ex:G-family}}{
  Let $X = \mathbb{Z}_{3}$ and let $G = S_{3} \cong \left\langle a, x \mid a^{3} = x^{2} = e,\ xax = a^{2}\right\rangle$.
  For each $g \in G$, $\ast^{g} : X \times X \to X$ is the binary operation defined by
  $$\begin{cases}
    a \ast^{g} b = 2b - a & (g = x,\ ax,\ ax^{2}),\\
    a \ast^{g} b = a & (\mbox{otherwise})
  \end{cases}$$
  for any $a, b \in X$ and $g \in G$.
  Then $\left(X, \left\{\ast^{g}\right\}_{g \in G} \right)$ is a $G$-family of racks. 
}

\prop[\cite{Ishii2015',Ishii-Matsuzaki-Murao2020}]{\label{Prop:Ass.MGR}}{
  Let $\left(X, \left\{\ast^{g}\right\}_{g \in G} \right)$ be a $G$-family of racks.
  Then $X \times G = \bigsqcup_{x \in X}\left(\left\{x\right\} \times G \right)$ is a multiple group rack with
  \begin{gather*}
    (x, g) \ast (y, h) = (x \ast^{h} y, h^{-1} gh),\ (x, g)(x, h) = (x, gh)
  \end{gather*}
  for any $x, y \in X$ and $g, h \in G$.
}\upshape

The multiple group rack $X \times G$ in Proposition~\ref{Prop:Ass.MGR} is called the \textit{associated multiple group rack} of the $G$-family of racks $\left(X, \left\{\ast^{g}\right\}_{g \in G} \right)$.

\example{\label{Ex:ass.MGR}}{
  Let $X$ be a rack.
  Define $n = \min \left\{k \in \mathbb{Z}_{>0} \mid \mbox{for any}\  x, y \in R, S_{y}^{k}(x) = x\right\}$.
  If the minimum integer does not exist, we set $n = \infty$.
  Then,
  $\left(X, \left\{\ast^{k}\right\}_{k \in \mathbb{Z}_{n}} \right)$ is a $\mathbb{Z}_{n}$-family of racks with operations $\left\{\ast^{k}\right\}_{k \in \mathbb{Z}_{n}}$ defined by $x \ast^{k} y = S_{y}^{k}(x)$ for any $x, y \in X$. 
  Here, when $n = \infty$, we regard $\mathbb{Z}_{\infty}$ as $\mathbb{Z}$.
  Thus,
  by   Proposition~\ref{Prop:Ass.MGR},
  $\textcolor{black}{X} \times \mathbb{Z}_{n} = \bigsqcup_{x \in X} (\left\{x\right\} \times \mathbb{Z}_{n})$ is the multiple group rack with
  \begin{equation*}
    (x, i) \ast (y, j) = (x \ast^{j} y, i),\quad (x, i)(x, j) = (x, i+j)
  \end{equation*}
  for any $x, y \in X$ and $i, j \in \mathbb{Z}_{n}$.
}

\section{\textcolor{black}{Coloring diagrams by multiple group racks}}{\label{Sec:MGRCols}}

Let $D$ be a Y-oriented diagram of a spatial surface.
An \textit{arc} of $D$ is a simple curve or a simple loop on $D$, obtained by cutting $D$ at undercrossings and vertices.
We denote the set of all arcs of $D$ by $\mathcal{A}(D)$.
\defn{\label{Def:MGRCols}}{
Let $X = \bigsqcup_{\lambda \in \Lambda} G_{\lambda}$ be a multiple group rack. 
An \textit{$X$-coloring} of a Y-oriented diagram $D$ or a \textit{coloring} of $D$ by $X$ is a map $C: \mathcal{A}(D) \to X$ satisfying the following conditions.
\begin{itemize}
  \item For each crossing of $D$, $C$ satisfies $C(a_{i}) \ast C(a_{j}) = C(a_{k})$,
  where $a_{i}, a_{j}, a_{k} \in \mathcal{A}(D)$ are as shown in the left side of Fig.~\ref{Fig:X-coloring_conditions}.
  \item For each vertex of $D$, $C$ satisfies $C(a_{i}), C(a_{j}), C(a_{k}) \in G_{\lambda}$ for some $\lambda \in \Lambda$ and
  $C(a_{i}) C(a_{j}) = C(a_{k})$, where $a_{i}, a_{j}, a_{k} \in \mathcal{A}(D)$ are as shown in the center or the right side of Fig.~\ref{Fig:X-coloring_conditions}.  
\end{itemize}
\begin{figure}[h]
  \centering
    \includegraphics{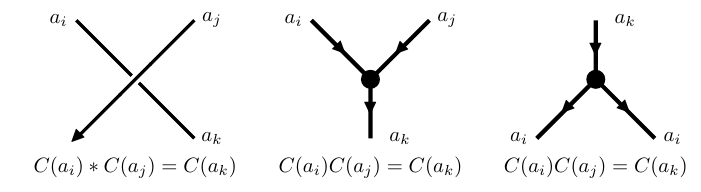}
    \caption{$X$-coloring conditions}{\label{Fig:X-coloring_conditions}}  
\end{figure}
We write $\mathrm{Col}_{X}(D)$ for the set of all $X$-colorings of $D$.
}

\thm[\cite{Ishii-Matsuzaki-Murao2020}]{\label{Thm:Colinv}}{
  Let $X$ be a multiple group rack and let $D$ and $D'$ be Y-oriented diagrams of spatial surfaces.
  If $F(D) \cong F(D')$,
  then there exists a bijection $\mathrm{Col}_{X}(D) \to \mathrm{Col}_{X}(D')$.
  In particular,
  the cardinality $\left\lvert \mathrm{Col}_{X}(D)\right\rvert$ is an invariant of the spatial surface $F(D)$.
}\upshape

We remark that $\left\lvert \mathrm{Col}_{X}(D)\right\rvert$ in Theorem~\ref{Thm:Colinv} does not depend on the choice of a Y-orientation of $D$.

\begin{figure}[h]
  \centering
\includegraphics{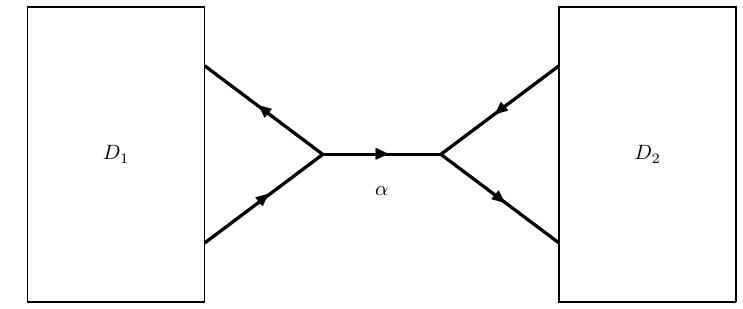}
\caption{$D = D_{1} \sharp_{\alpha} D_{2}$}{\label{Fig:connected_sum}}
\end{figure}

Let $D_{1}$ and $D_{2}$ be Y-oriented diagrams of spatial surfaces.
We denote by $D_{1} \sharp_{\alpha} D_{2}$ the Y-oriented diagram obtained from $D_{1}$ and $D_{2}$ by attaching an arc $\alpha$, as shown in Fig.~\ref{Fig:connected_sum}.
\question[cf.~\cite{Taniguchi2024}]{\label{Question:}}{
  Does there exist a multiple group rack $X = \bigsqcup_{\lambda \in \Lambda} G_{\lambda}$ which satisfies the following condition ($\star$)?
  \begin{itemize}
    \item[($\star$)] There exist a spatial surface $F$, \textcolor{black}{a} Y-oriented diagram $D = D_{1} \sharp_{\alpha} D_{2}$ of $F$, as shown in Fig.~\ref{Fig:connected_sum}, and \textcolor{black}{an} $X$-coloring $C$ of $D$
    such that $C(\alpha) \not\in \left\{e_{\lambda} \in G_{\lambda} \mid \lambda \in \Lambda\right\}$, 
    where $e_{\lambda}$ is the identity element of $G_{\lambda}$ ($\lambda \in \Lambda$).
  \end{itemize}
}\upshape

From Proposition~\ref{Prop:Ass.MGRCols},
\textcolor{black}{the} associated multiple group rack of \textcolor{black}{a} $G$-family of racks
does not satisfy the property $(\star)$ in Question~\ref{Question:}.

\prop{\label{Prop:Ass.MGRCols}}{
  Let $X \times G$ be the associated multiple group rack of a $G$-family of racks $\left(X, \left\{\ast^{g}\right\}_{g \in G} \right)$,
  and let $D = D_{1} \sharp_{\alpha} D_{2}$ be a Y-oriented diagram, as depicted in Fig.~\ref{Fig:connected_sum}. 
  Then, for any $X \times G$-coloring $C: \mathcal{A}(D) \to X \times G$,
  there exists an element $x \in X$ such that $C(\alpha) = (x, e)$,
  where $e$ is the identity element of $G$.
}\upshape

\begin{proof}
  Let $\phi: \mathcal{A}(D) \to G(N(F(D)))$ be the map that
  assigns to each arc of $D$ the corresponding meridian loop,
  where $N(F(D))$ is a regular neighborhood of the spatial surface $F(D)$, and 
  $G(N(F(D)))$ is the fundamental group of the exterior of $N(F(D))$.
  Then, there exists a group homomorphism $h: G(N(F(D))) \to G$ such that
  $\mathrm{pr} \circ C = h \circ \phi$, where $\mathrm{pr}: X \times G \to G$ is the projection defined by $\mathrm{pr}((x, g)) = g$ for any $(x, g) \in X \times G$.
  Since $\phi (\alpha)$ is the identity element of $G(N(F(D)))$ and $h$ is a group homomorphism,
  we have $\mathrm{pr} \circ C(\alpha) = h \circ \phi(\alpha) = e$.
  Hence, $C(\alpha) = (x, e)$ for some $x \in X$.
\end{proof}

\defn[\cite{Carter-Ishii-Saito-Tanaka2017,Matsuzaki-Murao2023}]{\label{Def:Ab_ext}}{
  Let $X = \bigsqcup_{\lambda \in \Lambda} G_{\lambda}$ be a multiple group rack (resp. multiple conjugation quandle).
  We define a free abelian group $C_{2}(X)$ by
  $$C_{2}(X) = \mathbb{Z} \left[ X^{2} \right] \sqcup \mathbb{Z} \left[\bigsqcup_{\lambda \in \Lambda} (G_{\lambda}^{2})\right].$$
  We present an element $(x, y) \in X^{2} \subset C_{2}(X)$ by $\left\langle x \right\rangle \left\langle y\right\rangle$ and 
  an element $(x_{1}, x_{2}) \in G_{\lambda}^{2} \subset C_{2}(X)$ by $\left\langle x_{1}, x_{2}\right\rangle$.
  Let $A$ be an abelian group.
  A homomorphism $f : C_{2}(X) \to A$ is a \textit{multiple group rack $2$-cocycle} (resp. \textit{multiple conjugation quandle $2$-cocycle}) if
  $X \times A = \bigsqcup_{\lambda \in \Lambda} (G_{\lambda} \times A)$ is a multiple group rack (resp. multiple conjugation quandle) with the following operations
  \begin{gather*}
    (x, a) \ast (y, b) = (x \ast y, a + f(\left\langle x\right\rangle \left\langle y\right\rangle)),\ 
    (x_{1}, a) (x_{2}, b) = (x_{1} x_{2}, a + b + f(\left\langle x_{1}, x_{2}\right\rangle))
  \end{gather*}
  for any $\lambda \in \Lambda$, $x, y \in X$, $x_{1}, x_{2} \in G_{\lambda}$, and $a, b \in A$.
  
  \textcolor{black}{Let $X$ be a multiple group rack, $A$ an abelian group, and $f: C_{2}(X) \to A$ a multiple group rack $2$-cocycle.}
  Then the multiple group rack (resp. multiple conjugation quandle) $X \times A$ is called an \textit{abelian extension} of $X$ by $A$.
}

\textcolor{black}{Let $X = \bigsqcup_{\lambda \in \Lambda} G_{\lambda}$ and $Y = \bigsqcup_{\mu \in M} H_{\mu}$ be multiple group racks.
A \textit{multiple group rack homomorphism} $f : X \to Y$ is a map satisfying the following conditions:}
\begin{itemize}
  \item[(i)] \textcolor{black}{For any $x, y \in X$, $f(x \ast y) = f(x) \ast f(y)$.}
  \item[(ii)] \textcolor{black}{For any $\lambda \in \Lambda$ and $a, b \in G_{\lambda}$, there exists $\mu \in M$ such that $f(a), f(b) \in H_{\mu}$ and $f(ab) = f(a) f(b)$.}
\end{itemize}

\lem{\label{Prop:Ab_Ext}}{
  \textcolor{black}{Let $X \times A$ be an abelian extension of a multiple group rack $X = \bigsqcup_{\lambda \Lambda} G_{\lambda}$ by an abelian group $A$ with respect to a multiple group rack $2$-cocycle $f: C_{2}(X) \to A$.
  Then the map $p: X \times A \to X,$ defined by $p((x, a)) = x$ for any $x \in X$ and $a \in A$, 
  is a multiple group rack homomorphism. 
}}\upshape

\begin{proof}
  \textcolor{black}{(i) For any $(x, a), (y, b) \in X \times A$,}
  \begin{equation*}
    \textcolor{black}{p((x, a) \ast (y, b)) = p((x \ast y, a + f(\left\langle x\right\rangle \left\langle y\right\rangle ))) = x \ast y\ \mbox{and}\ p((x, a)) \ast p((y, b)) = x \ast y.}
  \end{equation*}
  \textcolor{black}{Thus, we have $p((x, a) \ast (y, b)) = p((x, a)) \ast p((y, b))$.}

  \textcolor{black}{(ii) For any $\lambda \in \Lambda$, $x_{1}, x_{2} \in G_{\lambda}$, and $a, b \in A$,}
  \begin{equation*}
    \textcolor{black}{p((x_{1}, a)(x_{2}, b)) = p(x_{1} x_{2}, a + b + f(\left\langle x_{1}, x_{2}\right\rangle)) = x_{1} x_{2}\ \mbox{and}\ p((x_{1}, a)) p((x_{2}, b)) = x_{1} x_{2}.}
  \end{equation*}
  \textcolor{black}{Thus, we have $p((x_{1}, a)(x_{2}, b)) = p((x_{1}, a)) p((x_{2}, b))$.}

  \textcolor{black}{Therefore, it follows that $p: X \times A \to X$ is a multiple group rack homomorphism.}
\end{proof}


In \cite{Taniguchi2024},
it was shown that abelian extensions of multiple conjugation quandles, obtained from the knot quandles (resp. the knot $n$-quandles) of nontrivial knots by the method in Example~\ref{Ex:ass.MGR}, 
by $\mathbb{Z}$ (resp. $\mathbb{Z}_{n}$) 
satisfy the property $(\star)$.
For further details, 
see \cite{Carter-Ishii-Saito-Tanaka2017, Eisermann2003, Joyce1982, Matsuzaki-Murao2023, Matveev1982, Tanaka-Taniguchi2023}.

\thm[\cite{Taniguchi2024}]{\label{Thm:Taniguchi}}{
  There are infinitely many multiple conjugation quandles
  which satisfy the property ($\star$) in Question~\ref{Question:}.
}\upshape

\section{Main results}{\label{Sec:Main_results}}

The multiple group racks in Theorem~\ref{Thm:Taniguchi} are multiple conjugation quandles.
According to \cite{Ishii2015'}, 
multiple conjugation quandles are suitable for constructing invariants of handlebody-knots, handlebodies embedded in $S^{3}$.
However, they are not effective for distinguishing spatial surfaces.
Therefore,
for the purpose of constructing invariants of spatial surfaces, 
we require multiple group racks that are not multiple conjugation quandles.
In this paper,
we present a construction for multiple group racks that satisfy the property $(\star)$ in Question~\ref{Question:}, but are not multiple conjugation quandles.

\thm{\label{Thm:Main_result}}{
  Let $\left(X, \left\{\ast^{g}\right\}_{g \in G} \right)$ be a $G$-family of racks 
  and let $N \triangleleft G$ be a normal subgroup of $G$.
  Then $X \times (G \ltimes N) = \bigsqcup_{x \in X}\left(\left\{x\right\} \times \left(G \ltimes N\right) \right)$ is a multiple group rack defined by the following operations
  \begin{gather*}
    (x, (g_{1}, n_{1})) \ast (y, (g_{2}, n_{2})) = (x \ast^{g_{2} n_{2}} y, (g_{1}^{g_{2} n_{2}}, n_{1}^{g_{2} n_{2}})),\\
    (x, (g_{1}, n_{1})) (x, (g_{2}, n_{2})) = (x, (g_{1} g_{2}, n_{1}^{g_{2}} n_{2}))
  \end{gather*}
  for any $x, y \in X$ and $(g_{1}, n_{1}), (g_{2}, n_{2}) \in G \ltimes N$.
  Here, we use $g^{h}$ to denote the right conjugation action, i.e., $g^{h} = h^{-1}gh$ for any $g, h \in G$.
}\upshape

\begin{figure}[h]
  \centering
  \includegraphics{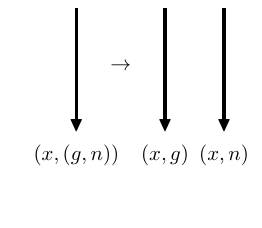}
  \includegraphics{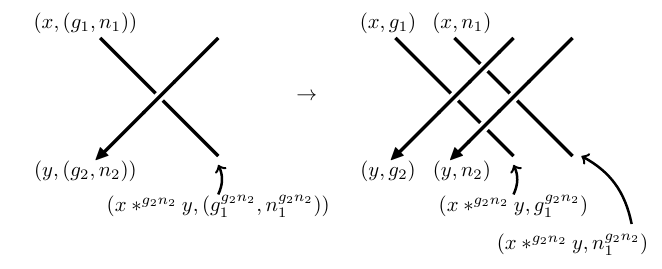}

  \includegraphics{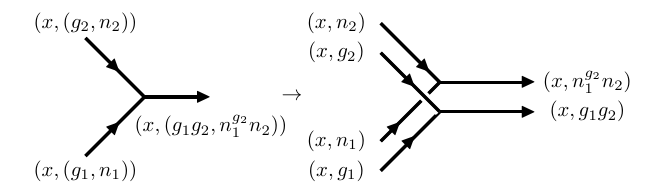}
  \caption{$x, y, x_{1}, x_{2} \in X, g, g_{1}, g_{2} \in G$, and $n, n_{1}, n_{2} \in N$}{\label{Fig:geometric meaning the two operations}}
\end{figure}

\rem{
  The operations in Theorem~\ref{Thm:Main_result} are motivated by the following diagrammatic interpretation illustrated in Fig.~\ref{Fig:geometric meaning the two operations}.
}

\begin{proof}
  We verify that $(X \times (G \ltimes N), \ast)$ satisfies the conditions (i)--(iii) of Definition~\ref{Def:MGR}.

  (i) For any $x, y \in X$ and $(g, n), (g_{1}, n_{1}), (g_{2}, n_{2}) \in G \ltimes N$,
  \begin{align*}
    (x, (g, n)) \ast \left( (y, (g_{1}, n_{1})) (y, (g_{2}, n_{2}))\right) &= (x, (g, n)) \ast (y, (g_{1} g_{2}, n_{1}^{g_{2}} n_{2}))\\
    &= (x \ast^{g_{1} g_{2} n_{1}^{g_{2}} n_{2}} y, g^{g_{1} g_{2} n_{1}^{g_{2}} n_{2}}, n^{g_{1} g_{2} n_{1}^{g_{2}} n_{2}})\\
    &= (x \ast^{g_{1} n_{1} g_{2} n_{2}} y, g^{g_{1} n_{1} g_{2} n_{2}}, n^{g_{1} n_{1} g_{2} n_{2}})\\
    &= (x \ast^{g_{1} n_{1}} y, g^{g_{1} n_{1}}, n^{g_{1} n_{1}}) \ast (y, (g_{2}, n_{2}))\\
    &= ((x, (g, n)) \ast (y, (g_{1}, n_{1}))) \ast (y, (g_{2}, n_{2})).
  \end{align*}
  For any $x, y \in X$ and $(g, n) \in G \ltimes N$,
  \begin{align*}
    (x, (g, n)) \ast (y, (e, e)) & = (x \ast^{e e} y, (g^{e e}, n^{e e}))\\
    &= (x, (g, n)),
  \end{align*}
  where $e$ is the identity element of $G$.

  (ii) For any $x, y, z \in X$ and $(g_{1}, n_{1}), (g_{2}, n_{2}), (g_{3}, n_{3}) \in G \ltimes N$,
  \begin{align*}
    ((x, (g_{1}, n_{1})) \ast (y, (g_{2}, n_{2}))) \ast (z, (g_{3}, n_{3})) &= (x \ast^{g_{2} n_{2}} y, (g_{1}^{g_{2} n_{2}}, n_{1}^{g_{2} n_{2}})) \ast (z, (g_{3}, n_{3}))\\
    &= ((x \ast^{g_{2} n_{2}} y) \ast^{g_{3} n_{3}} z, (g_{1}^{g_{2} n_{2} g_{3} n_{3}}, n_{1}^{g_{2} n_{2} g_{3} n_{3}})).
  \end{align*}
  \begin{align*}
    & ((x, (g_{1}, n_{1})) \ast (z, (g_{3}, n_{3}))) \ast ((y, (g_{2}, n_{2})) \ast (z, (g_{3}, n_{3}))) \\
    &= (x \ast^{g_{3} n_{3}} z, (g_{1}^{g_{3} n_{3}}, n_{1}^{g_{3} n_{3}})) \ast (y \ast^{g_{3} n_{3}} z, (g_{2}^{g_{3} n_{3}}, n_{2}^{g_{3} n_{3}}))\\
    &= ((x \ast^{g_{3} n_{3}} z) \ast^{(g_{2} n_{2})^{g_{3} n_{3}}} (y \ast^{g_{3} n_{3}} z), (g_{1}^{g_{3} n_{3} (g_{2} n_{2})^{g_{3} n_{3}}}, n_{1}^{g_{3} n_{3} (g_{2} n_{2})^{g_{3} n_{3}}}))\\
    &= ((x \ast^{g_{2} n_{2}} y) \ast^{g_{3} n_{3}} z, (g_{1}^{g_{2} n_{2} g_{3} n_{3}}, n_{1}^{g_{2} n_{2} g_{3} n_{3}})).
  \end{align*}
  Thus, we have $$((x, (g_{1}, n_{1})) \ast (y, (g_{2}, n_{2}))) \ast (z, (g_{3}, n_{3})) = ((x, (g_{1}, n_{1})) \ast (z, (g_{3}, n_{3}))) \ast ((y, (g_{2}, n_{2})) \ast (z, (g_{3}, n_{3}))).$$

  (iii) For any $x, y \in X$ and $(g, n), (g_{1}, n_{1}), (g_{2}, n_{2}) \in G \ltimes N$,
  \begin{align*}
    ((x, (g_{1}, n_{1})) (x, (g_{2}, n_{2}))) \ast (y, (g, n)) &= (x, (g_{1} g_{2}, n_{1}^{g_{2}} n_{2})) \ast (y, (g, n))\\
    &= (x \ast^{g n} y, ((g_{1} g_{2})^{gn}, (n_{1}^{g_{2}} n_{2})^{gn})).
  \end{align*}
  \begin{align*}
    ((x, (g_{1}, n_{1})) \ast (y, (g, n))) ((x, (g_{2}, n_{2})) \ast (y, (g, n))) &= (x \ast^{gn} y, (g_{1}^{gn}, n_{1}^{gn})) (x \ast^{gn} y, (g_{2}^{gn}, n_{2}^{gn}))\\
    &= (x \ast^{gn} y, (g_{1} g_{2})^{gn}, (n_{1}^{gn})^{g_{2}^{gn}} n_{2}^{gn})\\
    &= (x \ast^{gn} y, ((g_{1} g_{2})^{gn}, (n_{1}^{g_{2}} n_{2})^{gn})). \qedhere
  \end{align*}
\end{proof}

\thm{\label{Thm:Main_result2}}{
  In Theorem~\ref{Thm:Main_result},
  if the right conjugation action $N \curvearrowleft G$ is nontrivial,
  then the multiple group rack $X \times (G \ltimes N)$ satisfies the property ($\star$) in Question~\ref{Question:}.
}\upshape

\begin{proof}
  Assume that the right conjugation action $N \curvearrowleft G$ is nontrivial.
  Then, there exist elements $A \in G$ and $B \in N$ such that
  $\left[A, B\right] := ABA^{-1}B^{-1} \neq e$, 
  where $e$ is the identity element of $G$.
  Now, take the Y-oriented diagram $D$ and the coloring $C$ of $D$ by $X \times (G \ltimes N)$ as illustrated in Fig.~\ref{Fig:Main}.
  Then,
  $C(\alpha) = (x, (\left[A, B\right]^{-1}, \left[A, B\right])) \neq (x,(e, e)) \in \left\{x\right\} \times (G \ltimes N)$ for some $x \in X$.
  Hence,
  the multiple group rack $X \times (G \ltimes N)$ satisfies the property ($\star$) in Question~\ref{Question:}.
  \begin{figure}[h]
    \centering
    \includegraphics{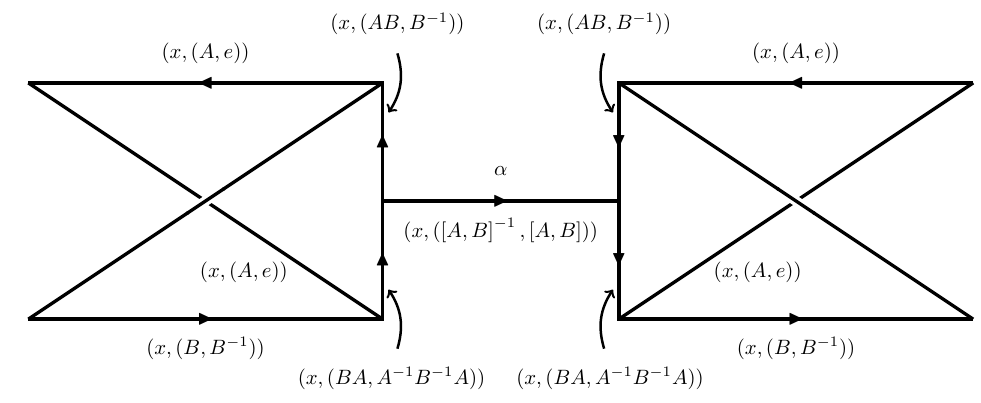}
    \caption{$D$ and $X \times (G \ltimes N$)-coloring $C$}{\label{Fig:Main}}
  \end{figure}\qedhere
\end{proof}

\prop{\label{Cor:}}{
  In Theorem~\ref{Thm:Main_result},
  if the right conjugation action $N \curvearrowleft G$ is nontrivial,
  we have the following statements (i) and (ii).
  \begin{itemize}
    \item[(i)] The multiple group rack $X \times (G \ltimes N)$ is not \textcolor{black}{the} associated multiple group rack of any $G$-family of racks.
    \item[(ii)] \textcolor{black}{The multiple group rack $X \times (G \ltimes N)$ is not an abelian extension of the associated multiple group rack $X \times G$ of the $G$-family of racks.}
  \end{itemize}
}\upshape

\begin{proof}
  (i) By Proposition~\ref{Prop:Ass.MGRCols},
  since the associated multiple group racks of $G$-families of racks do not satisfy the property $(\star)$,
  the multiple group rack $X \times (G \ltimes N)$, which satisfies the property $(\star)$, is not an associated multiple group rack of any $G$-family of racks.



\textcolor{black}{(ii) Let $p : X \times (G \ltimes N) \to X \times G$ be the map defined by $p(x, (g, n)) = (x, g)$.
  If the right conjugation action $N \curvearrowleft G$ is nontrivial,
  there exists $g \in G$ and $n \in N$ such that $g^{n} \neq g$.
  Then,
  for any $x \in X$,
  it follows that $p((x, (g, n)) \ast (x, (g, n))) = p(x \ast^{gn} x, (g^{n}, n^{gn})) = (x \ast^{gn} x, g^{n})$ and 
  $p((x, (g, n))) \ast p((x, (g, n))) = (x, g) \ast (x, g) = (x \ast^{g} x, g)$.
  Therefore,
  we have $p((x, (g, n)) \ast (x, (g, n))) \neq p((x, (g, n))) \ast p((x, (g, n)))$ for any $x \in X$.
  Hence $p$ is not a multiple group rack homomorphism.
  By Lemma~\ref{Prop:Ab_Ext},
  $X \times (G \ltimes N)$ is not an abelian extension of the multiple group rack $X \times G$ of the $G$-family of racks.}
\end{proof}

\textcolor{black}{In the following example,
we provide a pair of spatial surfaces that satisfies the following conditions:
(1) their boundaries are ambiently isotopic oriented links,
(2) their regular neighborhoods in $S^{3}$ are ambiently isotopic handlebody-knots,
(3) they cannot be distinguished by the number of colorings using the associated multiple group racks of any $G$-family of racks,
and (4) they can be distinguished by the number of colorings by a multiple group rack which is obtained by the construction (Theorem~\ref{Thm:Main_result}) in the present paper.}

\example{\label{Ex:}}{
  Let $D_{1}$ and $D_{2}$ be diagrams of spatial surfaces depicted in Fig.~\ref{Fig:ex}.
  Take the $G$-family of racks $\left(X, \left\{\ast^{g}\right\}_{g \in G} \right)$ as in Example~\ref{Ex:G-family}.
  Now, we set $N = G$ as the normal subgroup of $G$, 
  and consider the multiple group rack $X \times (G \ltimes G)$ 
  obtained by the construction in Theorem~\ref{Thm:Main_result}. 
  Assign Y-orientations to the diagrams $D_{1}$ and $D_{2}$ in Fig.~\ref{Fig:ex}. 
  Then, $\left\lvert\mathrm{Col}_{X \times (G \ltimes G)} (D_{1}) \right\rvert = 1458$ and $\left\lvert\mathrm{Col}_{X \times (G \ltimes G)} (D_{2}) \right\rvert = 1242$.
  Therefore, $F(D_{1}) \not\cong F(D_{2})$.
  \begin{figure}[h]
    \centering
    \includegraphics{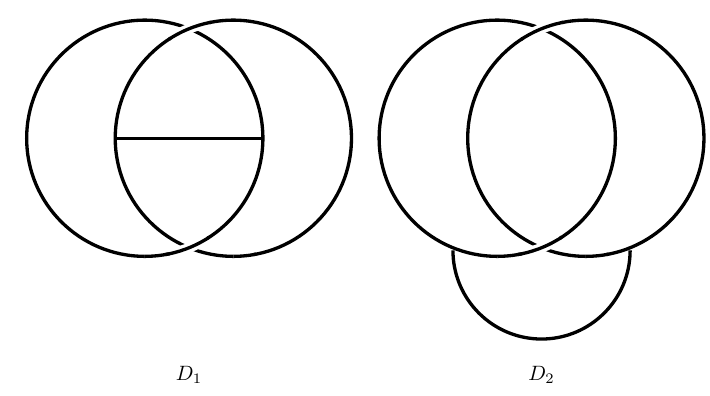}
    \caption{$D_{1}$ and $D_{2}$}{\label{Fig:ex}}
  \end{figure}
  On the other hand,
  some invariants are the same for the spatial surfaces $F(D_{1})$ and $F(D_{2})$ because
  $\partial F(D_{1})$ and $\partial F(D_{2})$ are ambiently isotopic to the link $\mathrm{L6n1}$,
  $N(F(D_{1}))$ and $N(F(D_{2}))$ are ambiently isotopic, and
  for any $G$-family of racks,
  the numbers of colorings $D_{1}$ and $D_{2}$ by the associated multiple group rack of the $G$-family of racks are the same.

  \textcolor{black}{Additionally, 
  by Proposition~\ref{Cor:},
  the multiple group rack $X \times (G \ltimes G)$ is neither the associated multiple group rack of any $G$-family of racks nor its abelian extension.
  Furthermore,
  the multiple group rack $X \times (G \ltimes G)$ is not a multiple conjugation quandle
  because $(0, (a, x)) \ast (0, (a, x)) = (0, (a^{2}, ax)) \neq (0, (a, x))$.}  
}

\section*{Acknowledgement}

The author would like to thank Seiichi Kamada and Yuta Taniguchi for discussions on this research.
This work was supported by JST SPRING, Grant Number JPMJSP2138.

\bibliographystyle{jplain}
\bibliography{grpd.bib}
\vspace{-3mm}

\address{(K. Arai) Department of Mathematics, Graduate School of Science, The University of Osaka, 1-1, Machikaneyama, Toyonaka, Osaka, 560-0043, Japan}

\email{u068111h@ecs.osaka-u.ac.jp}


\end{document}